\documentclass[12pt]{article}
\usepackage{latexsym}

\newcommand{\duk}{\noindent {\bf Proof. }}
\newcommand{\kduk}{\hfill $\Box$\bigskip}

\newcommand{\N}{\mathbf{N}}

\newcommand{\Q}{\mathbf{Q}}

\newtheorem{veta}{Theorem}[section]
\newtheorem{dusl}[veta]{Corollary}

\newtheorem{prop}[veta]{Proposition}

\def\cla#1#2#3#4#5#6{%autor, nazev, cas., rocnik, rok, stranky
  {\sc #1, }#2, {\it #3, }{\bf #4 }(#5), #6.}

\def\kni#1#2#3#4#5{%autor, nazev, nakladatel, sidlo, rok
  {\sc #1, }{\it #2, }#3, #4, #5.}

\begin{document}

\author{Martin Klazar\thanks{Department of Applied Mathematics (KAM) and Institute for Theoretical 
Computer Science (ITI), Charles University, Malostransk\'e n\'am\v est\'\i\ 25, 118 00 Praha, 
Czech Republic. ITI is supported by the project LN00A056 of the 
Ministry of Education of the Czech Republic. E-mail: {\tt klazar@kam.mff.cuni.cz}}}
\title{Counting set systems by weight}
\date{}

\maketitle
\begin{abstract}
Applying the enumeration of sparse set partitions, we show that the number of set systems 
$H\subset\exp(\{1,2,\dots,n\})$ 
such that $\emptyset\not\in H$, $\sum_{E\in H} |E|=n$ and 
$\bigcup_{E\in H} E=\{1,2,\dots,m\}$, $m\le n$, equals $(1/\log(2)+o(1))^nb_n$ where 
$b_n$ is the $n$-th Bell number. The same asymptotics holds if $H$ may be a multiset. 
If vertex degrees in $H$ are restricted to be at most $k$, the asymptotics is 
$(1/\alpha_k+o(1))^nb_n$ where $\alpha_k$ is the unique root of $\sum_{i=1}^k x^i/i!-1$ in $(0,1]$.
\end{abstract}

\section{Introduction}

Let $h_n'$, where $n\in\N=\{1,2,\dots\}$, be the number of simple set systems $H$ such that the edges 
$E\in H$ are nonempty finite subsets of $\N$, the {\em weight} $w(H)=\sum_{E\in H}|E|$ of $H$ equals 
$n$, and the vertex set is (normalized to be) $[m]=\{1,2,\dots,m\}$ 
for some $m$. Let $h_n''$ be the number of these set systems with weight $n$ if
edges may be repeated, i.e., if $H$ may be a multiset. For example, 
$H_1=\{\{1\},\{2\},\{3\}\}$,  $H_2=\{\{1,2\},\{3\}\}$, 
$H_3=\{\{1,3\},\{2\}\}$, $H_4=\{\{2,3\},\{1\}\}$, 
$H_5=\{\{1,2\},\{1\}\}$, $H_6=\{\{1,2\},\{2\}\}$, and $H_7=\{\{1,2,3\}\}$ show that $h_3'=7$. The three 
additional multisets $H_8=\{\{1\},\{1\},\{2\}\}$, $H_9=\{\{1\},\{2\},\{2\}\}$, and 
$H_{10}=\{\{1\},\{1\},\{1\}\}$ show that $h_3''=10$. 
We investigated $h_n'$ and $h_n''$ in \cite{klaz_dam} where we calculated the values
\begin{center}
\begin{tabular}{l||l|l|l|l|l|l|l|l|l|l}
$n$&1&2&3&4&5&6&7&8&9&10\\
\hline\hline
$h_n'$ & 1 & 2 & 7 & 28 & 134 & 729 & 4408 & 29256 & 210710 & 1633107\\
\hline
$h_n''$ & 1 & 3 & 10 & 41 & 192 & 1025 & 6087 & 39754 & 282241 & 2159916
\end{tabular}
\end{center}
and proved that $b_n\le h_n'\le h_n''\le 2^{n-1}b_n$, $b_n$ being the $n$-th 
Bell number (only the last inequality is nontrivial). In this note we prove 
the following stronger asymptotic bound. 

\bigskip\noindent
{\bf Theorem.} {\it If $n\to\infty$, $h_n'$ and $h_n''$ have the asymptotics
$$
((\log 2)^{-1}+o(1))^n\cdot b_n=(1.44269\dots+o(1))^n\cdot b_n
$$
where $b_n$ are Bell numbers.
}

\bigskip\noindent
The first Bell numbers are
\begin{center}
\begin{tabular}{l||l|l|l|l|l|l|l|l|l|l}
$n$&1&2&3&4&5&6&7&8&9&10\\
\hline\hline
$b_n$& 1 & 2 & 5 & 15 & 52 & 203 & 877 & 4140 & 21147 & 115975
\end{tabular}
\end{center}
and their asymptotics is recalled in Proposition~\ref{bell_asym}. 
In \cite{klaz_dam} we investigated besides $h_n'$ and $h_n''$ the numbers $h_n$ of simple set 
systems $H$ with vertex set $[n]$; $h_n$ were considered before by 
Comtet \cite[p. 165]{comt}, Hearne and Wagner 
\cite{hear_wagn}, and Macula \cite{macu} and they constitute entry 
A003465 of Sloane \cite{sloa_hand}. 

In Section 2 we prove the theorem. In Section 3 we give concluding comments, point out some analogies 
and pose some open questions. Now we recall and fix notation. For $n\in\N$, $\{1,2,\dots,n\}$ is 
denoted $[n]$. For $A,B\subset\N$, $A<B$ means that $x<y$ for every 
$x\in A$ and $y\in B$. We use notation $f(n)\ll g(n)$ as synonymous to the $f(n)=O(g(n))$ notation. The 
coefficient of $x^n$ in a power series $F(x)$ is denoted $[x^n]F$.  
A {\em set system} $H$ is here a finite multisubset of 
$\exp(\N)$ whose edges $E\in H$ are nonempty and finite. The {\em vertex set} 
is $V(H)=\bigcup_{E\in H}E$. The {\em degree} $\deg(x)=\deg_H(x)$ of a vertex $x\in V(H)$ in $H$ is 
the number of edges containing $x$.  If there are no multiple edges, 
we say that $H$ is {\em simple}. $H$ is a {\em partition} if its edges are mutually disjoint; 
in the case of partitions they are usually called {\em blocks}. The number of partitions $H$ 
with $V(H)=[n]$ is the Bell number $b_n$. The {\em weight} of a set system $H$ is 
$w(H)=\sum_{v\in V(H)}\deg(v)=\sum_{E\in H}|E|$. $H$ is {\em normalized} if $V(H)=[m]$ for some $m$.
In the proof of Proposition~\ref{simpleandall}
we work with more general set systems $H$ with vertex set contained in the dense linear 
order of fractions $\Q$. We {\em normalize} such set system by replacing it by the set system 
$H'=\{f(E):\ E\in H\}$, $V(H')=[m]$, where $f:V(H)\to[m]$ is the unique increasing bijection. 

\section{The proof}

To estimate $h_n'$ and $h_n''$ in terms of $b_n$, we transform a set system 
$H$ in a set partition with the same weight by splitting each vertex 
$v\in V(H)$ in $\deg_H(v)$ new vertices 
which are 1-1 distributed among the edges containing $v$. The following definitions and 
Propositions~\ref{ortrozkl} and \ref{tridyekv} make this idea precise. 

We call two set partitions $P$ and $Q$ of $[n]$ {\em orthogonal\/} if $|A\cap B|\le 1$ 
for every two blocks
$A\in P$ and $B\in Q$. $Q$ is an {\em interval\/} partition of $[n]$ if every block of $Q$ 
is a subinterval of $[n]$. For $n\in\N$ we define $W(n)$ to be the set of all pairs $(Q,P)$
such that $Q$ and $P$ are orthogonal set partitions of $[n]$ and $Q$ is moreover an interval partition. 
We define a binary relation $\sim$ on $W(n)$ by setting $(Q_1,P_1)\sim(Q_2,P_2)$ iff $Q_1=Q_2$ and there
is a bijection $f: P_1\to P_2$ such that for every $A\in P_1$ the blocks $A$ and $f(A)$ intersect the 
same intervals of the partition $Q_1=Q_2$. It is an equivalence relation. 

\begin{prop}\label{ortrozkl}
For every $n\in\N$, there is a bijection $(Q,P)\mapsto H(Q,P)$ between the set of equivalence classes 
$W(n)/\!\sim$ and the set $L(n)$ of normalized set systems $H$ with weight $n$. 
In particular, $h_n''=|L(n)|=|W(n)/\!\sim|$.
\end{prop}
\duk
We transform every $(Q,P)\in W(n)$, where $Q$ consists of the intervals 
$I_1<I_2<\dots<I_m$, in the set system 
$H=H(Q,P)=(E_A:\ A\in P)$ where $E_A=\{i\in [m]:\ A\cap I_i\ne\emptyset\}$. 
We have $w(H)=n$ and $V(H)=[m]$, so 
$H\in L(n)$. It is easy to see that equivalent pairs produce the same $H$ and nonequivalent pairs 
produce distinct 
elements of $L(n)$. 

Let $H\in L(n)$ with $V(H)=[a]$.
We split $[n]$ in $a$ intervals $I_1<I_2<\dots<I_a$ so that $|I_i|=\deg_H(i)$. 
For every $i\in[a]$ we fix arbitrary bijection $f_i: \{E\in H:\ i\in E\}\to I_i$. 
We define the partitions $Q=(I_1,I_2,\dots, I_a)$ and $P=(A_E:\ E\in H)$ where $A_E=\{f_i(E):\ i\in E\}$. 
Clearly, $(Q,P)\in W(n)$ and different choices of bijections $f_i$ lead to 
equivalent pairs. Also, $H(Q,P)=H$. Thus $(Q,P)\mapsto H(Q,P)$ is a bijection between 
$W(n)/\!\sim$ and $L(n)$. 
\kduk

The next proposition summarizes useful properties of the equivalence $\sim$ and the bijection 
$(Q,P)\mapsto H(Q,P)$. They follow in a straightforward way from the construction
and we omit the proof. 

\begin{prop}\label{tridyekv}
Let $(Q,P)\in W(n)$, $Q=(I_1<I_2<\dots<I_m)$, and $H=H(Q,P)$ (so $V(H)=[m]$). 
Then $\deg_H(i)=|I_i|$ for every
$i\in[m]$. The equivalence class containing $(Q,P)$ has at most $|I_1|!\cdot|I_2|!\cdot\dots\cdot|I_m|!$ 
pairs. It has exactly so many pairs if and only if $H$ is simple. 
\end{prop}

\begin{prop}\label{simpleandall}
For every $n\in\N$, $h_n'\le h_n''\le 2h_n'$.
\end{prop}
\duk
The first inequality is trivial. To prove the second inequality, we construct an injection from 
the set $N(n)$ of normalized non-simple set systems $H$ with weight $n$
in the set $M(n)$ of normalized simple set system 
$H$ with weight $n$. Then $h_n''=|M(n)|+|N(n)|\le 2|M(n)|=2h_n'$. 
We say that a vertex $v\in V(H)$ is {\em regular} if $\deg(v)\ge 2$ or if $v\in E$ 
for some $E\in H$ with $|E|\ge 2$, else we call $v$ {\em singular}. Thus $v$ is singular iff 
$\{v\}\in H$ and $\deg(v)=1$.  

Let $H\in N(n)$. We distinguish two cases. 
The first case is when every multiple edge of $H$ is a singleton. 
Then let $k\ge 2$ be the maximum multiplicity of an edge in $H$ and $v=u-1$ where 
$u\in V(H)$ is the smallest regular vertex in $H$; we may have $v=0$ and then $v$ is not a vertex 
of $H$.  We have  
$v<\max V(H)$ and insert between $v$ and $u$ new vertices $w_i$, $i=1,2,\dots,k-1$ and 
$v<w_1<w_2<\dots<w_{k-1}<u$. 
Then we replace every singleton multiedge $\{x\}$ with multiplicity $m$, $2\le m\le k$, 
(we have $x\ge u$) with the new single edge $\{w_1,w_2,\dots,w_{m-1},x\}$. 
Normalizing the resulting set system we get the set system $H'$. Clearly, $H'\in M(n)$. 

The second case is when at least one multiple edge in $H$ is not a singleton. We define $k$, $v$,
$u$, and $w_1,\dots, w_{k-1}$ as in the first case and replace every multiedge $E$ with multiplicity $m$, 
$2\le m\le k$, 
(we have $\min E\ge u$) by the new single edge $\{w_1,w_2,\dots,w_{m-1}\}\cup E$. 
We add between $w_{k-1}$ and $u$ a new vertex $s$ and add a new singleton edge $\{s\}$. 
This singleton edge is a marker discriminating between both cases and separating the new vertices 
$w_i$ from those in $E$. Since $m-1+|E|<m|E|$ if 
$|E|\ge 2$ and $m\ge 2$, the weight is still at most $n$. We add in the beginning sufficiently 
many new singleton edges
$\{-r\},\dots,\{-1\},\{0\}$ so that the resulting set system has weight exactly $n$. Normalizing it, we get 
the set system $H'$. Again, $H'\in M(n)$. Note that in both cases the least regular vertex in 
$H'$ is $w_1$ and that in both cases the longest interval in $V(H')$ that starts in $w_1$ and is 
a proper subset of an edge ends in $w_{k-1}$.

Given the image $H'\in M(n)$, in order to reconstruct $H$ we let $w\in V(H')$ 
be the least regular vertex (i.e., $w$ is the first vertex lying in an edge $E$ with $|E|\ge 2$) 
and let $I$ be the longest interval in $V(H')$ that starts in $w$ and is a proper subset of an edge. 
If $\max I+1$ is a singular vertex of $H'$, it must be $s$ and we are in the 
second case. Else there is no $s$ and we are in the first case. Knowing this and knowing (in the 
second case) which vertices are the dummy $w_i$, we uniquely reconstruct the multiedges of $H$. 
Thus $H\mapsto H'$ is an injection from $N(n)$ to $M(n)$. 
\kduk

For $k,n\in\N$ we define $h_{k,n}''$ to be the number of normalized set systems
with weight $n$ and maximum vertex degree at most $k$. The number of simple such systems is 
$h_{k,n}'$. Proposition~\ref{simpleandall_deg} can be proved by an injective argument similar to 
the previous one and we leave the proof to the 
interested reader as an exercise. But note that one cannot use the previous injection without change because 
it creates vertices with high degree.

\begin{prop}\label{simpleandall_deg}
For every $k,n\in\N$ we have $h_{k,n}'\le h_{k,n}''\le 2h_{k,n}'$.
\end{prop}

For the lower bound on $h_n''$ we need to count sparse partitions. A partition $P$ of $[n]$ 
is $m$-{\em sparse}, 
where $m\in\N$, if for every two elements $x<y$ of the 
same block we have $y-x\ge m$. Thus every partition is 1-sparse and 2-sparse partitions 
are those with no two consecutive numbers in the same block. If $m'<m$, 
every $m$-sparse partition is also 
$m'$-sparse. The number of $m$-sparse partitions of $[n]$ is denoted $b_{n,m}$. 
The following enumeration of sparse partitions was obtained by Prodinger 
\cite{prod} and Yang \cite{yang}, see also Stanley \cite[Problem 1.4.29]{stan}. Here we present a 
simple and nice proof due to Chen, Deng and Du \cite{chen_spol}.

\begin{prop}\label{sparse}
Let $m,n\in\N$. For $m>n$ there is only one $m$-sparse partition of $[n]$. 
For $m\le n$ the number $b_{n,m}$ of
$m$-sparse partitions of $[n]$ equals the Bell number $b_{n-m+1}$.
\end{prop}
\duk
For $m>n$ the only partition in question is that with singleton blocks. 
Let $P$ be a partition of $[n]$. We represent it by the graph $G=([n],E)$ where for $x<y$ we set 
$\{x,y\}\in E$ iff $x,y\in A$ for some block $A$ of $P$ and there is no $z\in A$ with $x<z<y$. 
The components of $G$ are increasing paths corresponding to the blocks of $P$. Equivalently, 
$G$ has the property 
that each vertex has degree at most $2$ and if it has degree $2$, it must lie between its two neighbors. 
Now assume that $P$ is 2-sparse. We transform 
$G$ in the graph $G'=([n-1],E')$ where $E'=\{\{x,y-1\}:\ \{x,y\}\in E,x<y\}$, i.e., we decrease the second
vertex of each edge by one. Note that $G'$ is again a graph (no loops arise). The property of $G$ is preserved by 
the transformation and hence the components of $G'$ are increasing paths and $G'$ describes a partition $P'$ of 
$[n-1]$. Clearly, $P$ is $m$-sparse iff $P'$ is $m-1$-sparse. Thus $P\mapsto P'$ maps $m$-sparse partitions of 
$[n]$ to $m-1$-sparse partitions of $[n-1]$. The inverse mapping is obtained simply by increasing the second 
vertex of each edge (in the graph describing a partition of $[n-1]$) by one. Thus $P\mapsto P'$ is a bijection 
between the mentioned sets. Iterating it, we obtain the stated identity. 
\kduk

\noindent
See \cite{chen_spol} for other applications of this bijection. We remark that the representing graphs
of partitions (but not the transformation of $G$ to $G'$) were used before by Biane \cite{bian} and
Simion a Ullman \cite{simi_ullm}. 

We need to compare, for fixed $m$, the growth of $b_n$ and $b_{n-m}$. The following asymptotics 
of Bell numbers is due to Moser and Wyman \cite{mose_wyma}.

\begin{prop}\label{bell_asym}
For $n\to\infty$, 
$$
b_n\sim\frac{\lambda(n)^{n+1/2}}{n^{1/2}\mathrm{e}^{n+1-\lambda(n)}}
$$ 
where the function $\lambda(n)$ is defined by $\lambda(n)\log\lambda(n)=n$.
\end{prop}

\noindent
It follows by a simple calculation that $b_{n-1}/b_n\sim\log n/n$. More generally, we have 
the following.

\begin{dusl}\label{podil}
If $m$ fixed and $n\to\infty$, 
$$
\frac{b_{n-m}}{b_n}\sim\bigg(\frac{\log n}{n}\bigg)^m.
$$
\end{dusl}

\noindent
In fact, a better approximation is $b_{n-1}/b_n\sim(\log n-\log\log n)/n$. Knuth 
\cite{knut} gives a nice account on the asymptotics of $b_n$ and shows that
$b_{n-1}/b_n=(\xi/n)(1+O(1/n))$ where $\xi\cdot\mathrm{e}^{\xi}=n$.  

We are ready to estimate the numbers of normalized set systems with weight 
$n$ and maximum degree at most $k$.

\begin{prop}\label{bound_deg}
For fixed $k\in\N$ and $n\to\infty$,
$$
\bigg(\frac{\log n}{n}\bigg)^{k-1} \bigg(\frac{1}{\alpha_k}\bigg)^n b_n\ll 
h_{k,n}''\ll
\bigg(\frac{1}{\alpha_k}\bigg)^n b_n
$$
where $\alpha_k$ is the only root of the polynomial $\sum_{i=1}^k x^i/i!-1$ in $(0,1]$.
\end{prop}
\duk
Let $i_{n,k}$ be the number of interval partitions $Q=(I_1<I_2<\dots<I_m)$ of $[n]$ such that 
$|I_i|\le k$ for all
$i$ and $Q$ is weighted by $(|I_1|!\cdot|I_2|!\cdot\dots\cdot|I_m|!)^{-1}$. It follows that
$$
i_{n,k}=[x^n]\frac{1}{1-\sum_{i=1}^k x^i/i!}\sim c_k\bigg(\frac{1}{\alpha_k}\bigg)^n 
$$
with some constant $c_k>0$ because $\alpha_k$ is the only root of the denominator in $(0,1]$ and it is simple.
Using Propositions \ref{ortrozkl}, \ref{tridyekv}, and \ref{simpleandall_deg}, we obtain the inequalities
$$
i_{n,k}b_{n,k}\le h_{k,n}''\le 2i_{n,k}b_n.
$$
In the first inequality we use the fact that if $Q$ is an interval partition of $[n]$ with interval 
lengths 
at most $k$ and $P$ is a $k$-sparse partition of $[n]$, then $Q$ and $P$ are always orthogonal. In the 
second inequality we neglect orthogonality of the pairs $(Q,P)$ but we assume that the 
corresponding equivalence classes in $W(n)$ have
full cardinalities $|I_1|!\cdot|I_2|!\cdot\dots\cdot|I_m|!$. By Proposition~\ref{tridyekv}, this gives an upper 
bound for $h_{k,n}'$. Using Proposition~\ref{simpleandall_deg}, we get an upper bound for $h_{k,n}''$.
The explicit lower and upper bounds on $h_{k,n}''$ now follow from the above asymptotics of $i_{n,k}$, 
Proposition~\ref{sparse}, and Corollary~\ref{podil}.
\kduk

Note that $1/\alpha_2=(1+\sqrt{3})/2$. Thus we have roughly 
$((1+\sqrt{3})/2)^nb_n=(1.36602\dots)^nb_n$ normalized set systems with weight $n$, in which each vertex lies 
in one or two edges. 

\begin{veta}
For $n\to\infty$,
$$
h_n''=\bigg(\frac{1}{\log 2}+o(1)\bigg)^nb_n=(1.44269\dots+o(1))^nb_n.
$$
\end{veta}
\duk
Let $i_n$ be the number of interval partitions $Q$ of $[n]$, weighted as in the previous proof. 
As in the case 
of bounded degree, by Propositions \ref{ortrozkl}, \ref{tridyekv}, and \ref{simpleandall} we have 
the upper bound
$$
h_n''\le 2i_nb_n\sim c\bigg(\frac{1}{\log 2}\bigg)^n b_n
$$
because 
$$
i_n=[x^n]\frac{1}{1-\sum_{i=1}^{\infty} x^i/i!}=[x^n]\frac{1}{2-\mathrm{e}^x}
$$
and $\log 2$ is a simple zero of $2-\mathrm{e}^x$. As for the lower bound, $h_n''\ge h_{k,n}''$ for every 
$k,n\in\N$. It is easy to show that $\alpha_k\downarrow\log 2$ for $k\to\infty$. Hence, by the lower bounds in 
Proposition~\ref{bound_deg}, for any $\varepsilon>0$ we have $h_n''>((\log 2)^{-1}-\varepsilon)^nb_n$ for 
$n$ big enough.
\kduk

\section{Concluding remarks}

P. Cameron investigates in \cite{came_draft} a family of enumerative problems on 0-1 matrices including 
$h_n'$ and $h_n''$ as particular cases. He defines $F_{ijkl}(n)$, $i,j,k,l\in\{0,1\}$, 
to be the number of 
rectangular 0-1 matrices with no zero row or column and with $n$ 1's, where $i=0$, resp. $i=1$, 
means that matrices differing only by a permutation of rows are identified, resp. are considered as 
different; $j=0$, resp. $j=1$, means that matrices 
with two equal rows are forbidden, resp. are allowed; and the values of $k,l$ refer to the same 
(non)restrictions for columns. Notice that $F_{ijkl}(n)$ is nondecreasing in each of the 
arguments $i,j,k,l$. Representing set systems by incidence matrices, rows standing for 
edges and columns for vertices, we see that $h_n''=F_{0111}(n)$ and $h_n'=F_{0011}(n)$. 
(Another language which can be used to deal with these problems, besides set systems and 0-1 matrices, 
is that of bipartite graphs.) In \cite{came_draft} it is shown that 
$F_{1111}(n)\sim Ac^{n+1}n!$ where $A=\frac{1}{4}\exp(-(\log 2)^2/2)\approx 0.19661$ and 
$c=(\log 2)^{-2}\approx 2.08137$. $F_{0101}(n)$ is A049311 of \cite{sloa_hand}, see also 
Cameron \cite{came}. P. Cameron asks in \cite[Problem 3]{came_pr} if there is an effective
algorithm to calculate $F_{0101}(n)$; for $h_n'$ and $h_n''$ such algorithms are given in 
\cite{klaz_dam}.

Interestingly, in the so far derived asymptotics of the functions $F_{ijkl}(n)$ the constant 
$\log 2\approx 0.69314$ appears quite often. Our theorem says that allowing intersecting  
blocks in ``partitions'' of $[n]$ is reflected in the counting function by multiplying 
$b_n$ by the exponential factor $(\log 2)^{-n}$. The same phenomenon occurs for 
counting injections and surjections. If $i_n$ is the number of injections from $[n]$ to 
$\N$ with images normalized to $[n]$ and $s_n$ is the number of all mappings (``injections with 
intersections'') from $[n]$ to $\N$, again with images normalized to $[m]$ (i.e., $s_n$ counts
surjections from $[n]$ to $[m]$), then $i_n=n!$ 
(trivial) and $s_n\sim c(\log 2)^{-n}n!$ where $c=(2\log 2)^{-1}$ 
(a nice exercise on exponential generating functions, see Flajolet and Sedgewick 
\cite[Chapter 2.3.1]{flaj_sedg}). 

Another parallel can be led between sparse partitions and sparse words. We say that a word 
$u=a_1a_2\dots a_l$ over an alphabet $A$, $|A|=r$, is $k$-sparse if $a_i=a_j$, $j>i$, implies $j-i\ge k$. 
(We remark that $k$-sparse words are basic objects in the theory of generalized Davenport--Schinzel 
sequences, see Klazar \cite{klaz_02}). The two notions of sparseness, in fact, coincide: 
$u=a_1a_2\dots a_l$ defines
a partition $P$ of $[l]$ via the equivalence $i\sim j\iff a_i=a_j$ and then, obviously, 
$u$ is $k$-sparse if and only if $P$ is $k$-sparse. A given partition of $[l]$ can be defined by 
many words $u$ (even if $A$ is fixed). The unique canonical defining words are {\em restricted growth 
strings}, see \cite{knut} for their properties and more references. (Another term for 2-sparse words 
is {\em Smirnov words} which is used, e.g., in \cite{flaj_sedg} or in Goulden and Jackson 
\cite{goul_jack}; we prefer the former term.) If $v_n$ is the number of all words over $A$
($|A|=r$) with length $n$ and $s_{k,n}$ is the number of those which are $k$-sparse, then 
$v_n=r^n$ (trivial) and
\begin{eqnarray*}
s_{k,n}&=&r(r-1)\dots (r-k+2)(r-k+1)^{n-k+1}\\
&=&\frac{r(r-1)\dots (r-k+2)}{(r-k+1)^{k-1}}\left(1-\frac{k-1}{r}\right)^nv_n
\end{eqnarray*} 
(by a simple direct counting, for the generating functions approach see \cite[Chapter 3.6.3]{flaj_sedg}). 
In the case of words over a fixed alphabet sparseness diminishes the counting function by an 
exponential factor. Fortunately, for partitions the decrease is much smaller, only by a polynomial factor 
(Corollary~\ref{podil}). 

We conclude by mentioning two natural questions. What is the precise asymptotics of 
$h_{k,n}''$ and $h_n''$?
By Propositions~\ref{simpleandall} and \ref{simpleandall_deg}, $1/2\le h_{k,n}'/h_{k,n}''\le 1$ and 
$1/2\le h_n'/h_n''\le 1$. Does this ratio go to $1$ as $n\to\infty$?

\end{document}